\documentclass[10pt,a4paper]{article}
\usepackage[utf8]{inputenc}
\usepackage[english]{babel}
\usepackage{amsmath}
\usepackage{amsfonts}
\usepackage{amssymb}
\usepackage{graphicx,color,float,verbatim}
\usepackage[affil-it]{authblk}

\definecolor{verde}{cmyk}{1,0,1,0.3}
\definecolor{naranja}{cmyk}{0,0.7,0.8,0}
\newcommand{\nj}{}

\newtheorem{theorem}{Theorem}

\newcommand{\R}{\mathbb{R}}

\begin{document}
\date{}
\small{
\title{Existence of solutions for a biological model using topological degree theory\\
{\large{Existencia de soluciones para un modelo biol\'ogico empeando la teor\'ia de grado topol\'ogico}}}
\author{Carlos H\'ector Daniel Alliera}
\affil{Departamento de Matem\'{a}tica, FCEyN - 
Universidad de Buenos Aires,
Pabell\'on I, Ciudad Universitaria, Buenos Aires, Argentina,
\textcolor{blue}{calliera@dm.uba.ar, http://cms.dm.uba.ar}}
\maketitle 
\begin{abstract}

\begin{small}
Topological degree theory is a useful tool for studying systems of differential equations. In this work, a biological model is considered. Specifically, we prove the existence of positive 
 $T$-periodic solutions of a system of delay differential equations 
for a model with feedback arising on Circadian oscillations in the Drosophila period gene protein. 

\begin{center}

\textbf{Resumen}
\end{center}
La teor\'ia de grado topológico es una herramienta \'util para estudiar sistemas de ecuaciones diferenciales. En este trabajo analizamos un modelo biol\'ogico; espec\'ificamente, 
probamos la existencia de soluciones positivas $T$-peri\'odicas de un sistema de ecuaciones diferenciales con retardo basado en el modelo auto-regulado de los ciclos Circadianos de prote\'inas a nivel gen\'etico de la mosca de la fruta (\textit{Drosophila}).
\\
\textit{Palabras Clave}: Ecuaciones diferenciales con retardo; Soluciones peri\'odicas, Modelos auto-regulados; Grado topol\'ogico.\\
\end{small} 
\end{abstract}

\section{Introduction}


Let us consider a model proposed by Goldbeter \cite{gB}, 
who showed the variation on PER: Period of \textit{messenger of Ribo-Nucleic Acid (mRNA)} in Drosophila (often called ``fruit flies") related to circadian rhythms.
{Here, a nonautonomous version of 
the model is considered with the aim of proving the existence of periodic solutions 
by means of a powerful topological tool: 
the Leray-Schauder degree.}
{In the original model, the existence of a positive steady state can be shown, under appropriate conditions, by the use of the Brouwer degree. 
As we shall see, when the parameters are replaced by 
periodic functions, 
essentially the same conditions yield the existence of positive periodic solutions.}  

\ \
\begin{figure}[h!]
\begin{center}
\includegraphics[scale=.5]{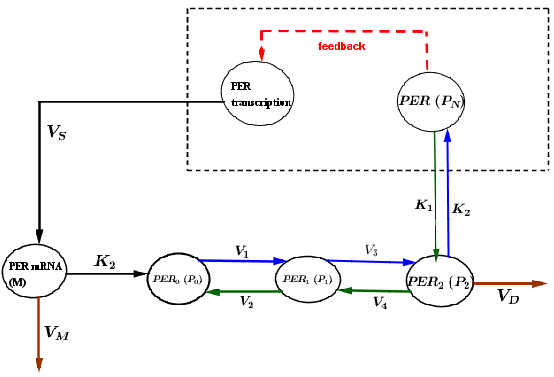}
\caption{{\small The model
for the circadian variation in PER.}}
\end{center}
\end{figure}

\section{The model}
\subsection{General features}
\begin{enumerate}
\item[\nj{I)}] This negative feedback will be described by an equation of Hill type
in which $n$ denotes the degree of cooperativity, and $K(t)$, the threshold
repression function.
\item[\nj{II)}]  To simplify the model, we consider that $P_N$ behaves
directly as a repressor.
\item[\nj{III)}] The constants $K,K_i$ and $V_j$ denote the maximum rate and
Michaelis constant of the kinase(s) and phosphatase(s) involved in the
reversible phosphorylation of $P_0$, into $P_1$, and of $P_1$, into $P_2$ are not negative.
\item[\nj{IV)}] Maximum accumulation rate of cytosol is denoted by $V_s$.
\item[\nj{V)}]  Cytosol is degraded enzymically,
in a Michaelian manner, at a maximum rate $V_m$.
\item[\nj{VI)}] Functions of this system are:
\begin{enumerate}
\item Cytosolic concentration is denoted by $M$.
\item We consider only three states of the protein:
unphosphorylated ($P_0$), monophosphorylated ($P_1$) and bisphosphorylated
($P_2$). 
\item Fully phosphorylated
form of PER ($P_2$) is degraded in a Michaelian manner, at a
maximum rate $V_d$ and also transported into the nucleus, at a rate characterized
by the apparent first-order rate constant  $k_1$.
\end{enumerate}
\item[\nj{VII)}] The rate of synthesis of PER, proportional to $M$, is characterized by
an apparent first-order rate constant $K_s$.
\item[\nj{VIII)}] Transport of the nuclear, bisphosphorylated form of PER ($P_N$) into the cytosol is characterized by the apparent first-order rate constant $k_2$.
\item[\nj{IX)}]  The model could be readily extended to include a larger number
of phosphorylated residues.
\end{enumerate}
{With this in mind, our non-autonomous version of Goldbeter's system reads:} 
\begin{equation}
\label{t1}
\begin{array}{l}
\dfrac{dM}{dt}=\dfrac{V_S(t)K_1(t)^{n}}{K_1^{n}(t)+P_N(t)^{n}}-\dfrac{V_{m}(t)M(t)}{K_{m_1}(t)+M(t)}, \\
\\

\dfrac{dP_0}{dt}= K_s(t) M(t)+\dfrac{V_2(t) P_{1}(t)}{K_2(t)+P_1(t)}-\dfrac{V_1(t)P_{0}(t)}{K_1(t)+P_0(t)}, \\
\\

\dfrac{dP_1}{dt}=\dfrac{V_1(t) P_{0}(t)}{K_1(t)+P_0(t)} + \dfrac{V_4(t)P_2(t)}{K_4(t)+P_2(t)}-P_{1}(t)\left( \dfrac{V_2(t)}{K_2(t)+P_1(t)}+\dfrac{V_3(t)}{K_3(t)+P_1(t)}\right), \\
\\

\dfrac{dP_2}{dt}= \dfrac{V_3(t)P_1(t)}{K_3(t)+P_1(t)}+k_{2}(t)P_N(t)-P_{2}(t)\left( k_1(t)+\dfrac{V_4(t)}{K_4(t)+P_2(t)}+\dfrac{V_d(t)}{K_d(t)+P_2(t)}\right), \\
\\

\dfrac{dP_N}{dt}= k_1(t) P_2(t)-k_{2}(t)P_N(t)
\end{array}
\end{equation}
where $K_i,\ i=1,2,3,4,d,m_1$, $k_1,\ k_2$ and $V_j,\ j=1,2,3,4,S,m,d$ are strictly positive, continuous $T$-periodic} functions.
We shall prove that, under accurate assumptions 
to be specified below, the system admits at least one 
positive $T$-periodic solution.

\section{Existence of positive periodic solutions}

In order to apply the topological degree method to problem 
(\ref{main}), 
let us consider the 
space of continuous $T$-periodic vector functions 
$$C_{T}:=\{u \in C(\mathbb{R},\mathbb{R}^{5}):u(t)=u(t+T)\,\hbox{ for all $t$}\},$$
equipped with the standard uniform norm,
and the positive cone 
$$\mathcal{K}:=\{u\in C_T: u_j \ge 0, j=1\ldots, 5\}.$$
Thus, the original problem can be written as 
$Lu=Nu$, where 
$L:C ^1\cap C_T\to C$ is given by $Lu:=u'$ and the nonlinear operator 
$N:\mathcal{K}\to C_T$ is 
defined as the right-hand side of system (\ref{t1}).
For convenience, the average of a function $u$ shall be denoted by $\overline u$, namely 
$\overline u:=\frac{1}{T}\int_0^T u(t)\, dt$. 
Also, identifying $\R^{5}$ with the subset of constant functions of $C_T$, we may define the function $\phi:[0,+ \infty)^{5}\to \R^{5}$ given by
$\phi(x) := \overline {Nx}$.

For the reader's convenience, let us 
summarize the basic properties of the Leray-Schauder degree which, roughly speaking, can be regarded as an algebraic count 
of the zeros of a mapping $F:\overline \Omega\to E$, where 
$E$ is a Banach space and $\Omega\subset E$ is open and bounded. 
In more precise terms, assume that $F=I-K$, where $K$ is compact 
and $F\neq 0$ on $\partial \Omega$. 
The degree $deg_{LS}(F,\Omega,0)$ is defined as the Brouwer degree $deg_B$ of its restriction 
$F|_V:\Omega \cap V\to V$, where $V$ is an accurate finite-dimensional subspace of $E$. 
In particular, if the range of $K$ is finite dimensional, then one may take $V$ as the subspace spanned by ${\rm Im}(K)$. 
If $deg_{LS}(F,\Omega,0)$ is different from $0$, 
then $F$ vanishes in $\Omega$; moreover, the degree is invariant over a continuous homotopy 
$F_\lambda:=I-K_\lambda$ with $K_\lambda$ compact and $F_\lambda\neq 0$ over $\partial \Omega$. 
Finally, we recall that 
if $T:\R^n\to\R^n$ is a diffeomorphism and $0\in T(A)$ for some open bounded $A\subset \R^n$ then 
$deg_B(T,A,0)$ is just the sign of the jacobian determinant of $T$ at the (unique) pre-image of $0$. 
The following continuation theorem is a direct consequence of the standard
topological degree methods (see e.g. \cite{PA}).

\begin{theorem}

$\label{TC}$
Assume there exists {$\Omega\subset \mathcal{K}^\circ$} open and bounded such that:

\begin{enumerate}
\item[a)] The problem $Lu=\lambda Nu$ has no solutions on $\partial\Omega$ for $0<\lambda<1$.
\item[b)] $\phi(u)\neq 0$ for all {$u\in\partial\Omega \cap\mathbb{R}^{5}$}.
\item[c)] {$deg_B(\phi,\Omega \cap\mathbb{R}^{5}, 0)\neq 0$}. 

\end{enumerate}
Then $(\ref{t1})$ has at least one solution in {$\overline \Omega$}.
\end{theorem}

\subsection{A priori bounds}

In this section, we shall find appropriate bounds for 
the solution of the problem $Lu=\lambda Nu$ with $\lambda\in (0,1)$. For convenience, let us fix the following notation 
for the 
minima and maxima 
of all the functions involved in the model,namely
$$0<\mathit{v}_i\leq V_i(t)\leq\mathcal{V}_i,\ 0<\kappa_j\leq K_j(t)\leq\mathcal{K}_j,\ 0<\hat{k}_l\leq k_l(t)\leq \Bbbk_l ,\ \forall\ i,\ j,\ l.$$
Now assume that $u\in \mathcal K^\circ$ satisfies  $Lu=\lambda Nu$ for some $0<\lambda<1$.\\
Let us firstly consider a value $t^*$ where $M$ achieves an absolute maximum, then $M'(t^*)=0$ and hence
$$\dfrac{V_S(t^*)K_1(t^*)^{n}}{K_1^{n}(t^*)+P_N(t^*)^{n}}=\dfrac{V_{m}(t^*)M(t^*)}{K_{m_1}(t^*)+M(t^*)}
\ge\dfrac{\mathit{v}_{m}M(t^*)}{\mathcal{K}_{m_1}+M(t^*)}
:=b_M(M(t^*)).$$
If 
\begin{equation} \label{hip1}
\mathit{v}_{m}> \mathcal V_S
\end{equation} 
then 
$$M(t^*)=b_M^{-1}\left( \dfrac{V_S(t^*)K_1(t^*)^{n}}{K_1^{n}(t^*)+P_N(t^*)^{n}}\right)<b_M^{-1}\left( V_S(t^*)\right)\leq \frac{\mathcal{V}_S\mathcal{K}_{m_1}}{\mathit{v}_M-\mathcal{V}_S}:=\mathcal{M} $$
Next, suppose that $P_0$ achieves its absolute maximum at some point, denoted again $t^*$, then 
$$K_s(t^*) M(t^*)+\dfrac{V_2(t^*) P_{1}(t^*)}{K_2(t^*)+P_1(t^*)}=\dfrac{V_1(t^*)P_{0}(t^*)}{K_1(t^*)+P_0(t^*)}
\ge 
\dfrac{\mathit{v}_1P_{0}(t^*)}{\mathcal K_1+P_0(t^*)}:=b_0(P_0(t^*)).$$
Thus, under the condition 
\begin{equation} \label{hip2}
\mathcal{K}_S\mathcal{M} + \mathcal{V}_2 <v_1,
\end{equation}
we deduce that 
$$P_0(t^*)=b_0^{-1}\left(K_s(t^*) M(t^*)+\dfrac{V_2(t^*) P_{1}(t^*)}{K_2(t^*)+P_1(t^*)} \right) < \frac {\mathcal{K}_S\mathcal{M} + \mathcal{V}_2}{v_1 -
(\mathcal{K}_S\mathcal{M} + \mathcal{V}_2)}\mathcal{K}_1:=\mathcal{P}_0.$$
Next, an upper bound $\mathcal P_1$ for $P_1$ is readily obtained 
in the following way. Let us denote again by $t^*$ 
a value where $P_1$ achieves its absolute maximum, then
$$\dfrac{V_1(t^*) P_{0}(t^*)}{K_1(t^*)+P_0(t^*)} + \dfrac{V_4(t^*)P_2(t^*)}{K_4(t^*)+P_2(t^*)}=P_{1}(t^*)\left( \dfrac{V_2(t^*)}{K_2(t^*)+P_1(t^*)}+\dfrac{V_3(t^*)}{K_3(t^*)+P_1(t^*)}\right).$$
When 
$P_1(t^*)\gg 0$, the right-hand side gets close to 
$V_2(t^*) + V_3(t^*)$, while the left-hand side is always less or equal than $\frac{\mathcal{V}_1\mathcal{P}_0}{\kappa_1+\mathcal{P}_0} + \mathcal{V}_4$. 
Thus, the existence 
of $\mathcal P_1$ is guaranteed by the condition
\begin{equation}
\label{hip3}
\frac{\mathcal{V}_1\mathcal{P}_0}{\kappa_1+\mathcal{P}_0} + \mathcal{V}_4 < \min_{t\in \mathbb R}\lbrace V_2(t)+V_3(t)\rbrace.
\end{equation}
The remaining upper bounds are obtained as follows. 
In the first place, define a 
new variable $Q:=P_N+P_2$ which satisfies the equation:
$$\dfrac{dQ}{dt}=\dfrac{V_3(t)P_1(t)}{K_3(t)+P_1(t)}-P_{2}(t)\left(\dfrac{V_4(t)}{K_4(t)+P_2(t)}+\dfrac{V_d(t)}{K_d(t)+P_2(t)}\right).$$
If $Q$ achieves its absolute maximum at $t^*$, then
$$\dfrac{\mathcal V_3\mathcal P_1}{\kappa_3+\mathcal P_1} \ge 
\dfrac{V_3(t^*)P_1(t^*)}{K_3(t^*)+P_1(t^*)}- P_{2}(t^*)\left(\dfrac{V_4(t^*)}{K_4(t^*)+P_2(t^*)}+\dfrac{V_d(t^*)}{K_d(t^*)+P_2(t^*)}\right).$$
As before, if the condition
\begin{equation}
\label{hip4}
\dfrac{\mathcal V_3\mathcal P_1}{\kappa_3+\mathcal P_1}<\min_{t\in\mathbb R}(V_4(t)+V_d(t))
\end{equation} 
is assumed, then 
$P_2(t^*)\le \tilde P$ for some $\tilde P$. 
Moreover, from the fourth equation of the system we deduce the existence of a constant $C$ such that 
$\dfrac{dP_2}{dt}\ge -CP_2(t)$. 
Hence we obtain, for all $t$, that  
$P_2(t)\le e^{CT}\tilde P:=\mathcal P_2$. This provides also an upper bound for $Q(t)$ and, consequently, 
an upper bound $\mathcal P_N$ for $P_N(t)$. 

After upper bounds are established, we proceed with the lower bounds as follows. 
Assume that $M$ achieves its absolute minimum at some
$t_*$, then we use again the fact that $M'(t_*)=0$ to obtain:
 $$\dfrac{V_{m}(t_*)M(t_*)}{K_{m_1}(t_*)+M(t_*)}=\dfrac{V_S(t_*)K_1(t_*)^{n}}{K_1^{n}(t_*)+P_N(t_*)^{n}}\ge 
\dfrac{\mathit{v}_S\kappa_1^{n}}{\kappa_1^{n}+\mathcal P_N^{n}}.$$
This shows that 
$M_1(t)\ge \mathfrak{m}$ for some positive constant 
$\mathfrak{m}$. In the same way, we find a lower bound $\mathfrak{p}_0$ for
$P_0$ using the fact that
$$\dfrac{V_1(t_*)P_{0}(t_*)}{K_1(t_*)+P_0(t_*)}=
K_s(t_*) M(t_*)+\dfrac{V_2(t_*) P_{1}(t_*)}{K_2(t_*)+P_1(t_*)}\ge 
\kappa_s\mathfrak m.$$
Next, suppose that 
$P_1$ achieves its absolute minimum at $t_*$, then
$$P_{1}(t_*)\left( \dfrac{V_2(t_*)}{K_2(t_*)+P_1(t_*)}+\dfrac{V_3(t_*)}{K_3(t_*)+P_1(t_*)}\right)
>\dfrac{V_1(t_*) P_{0}(t^*)}{K_1(t_*)+P_0(t_*)}\ge 
\dfrac{v_1 \mathfrak{p}_{0}}{\mathcal K_1+\mathfrak{p}_0}>0
$$
which yields the existence of a positive lower bound $\mathfrak p_1$.
Finally, positive lower bounds for $P_2$ and $P_N$ are obtained by means of the function $Q=P_2+P_N$. Indeed, if $Q$ achieves its absolute minimum at some $t_*$, then 
$$
P_{2}(t_*)\left(\dfrac{V_4(t_*)}{K_4(t_*)+P_2(t_*)}+\dfrac{V_d(t_*)}{K_d(t_*)+P_2(t_*)}\right)\ge 
\dfrac{v_3 \mathfrak{p}_{1}}{\mathcal K_3+\mathfrak{p}_1}
$$
and we deduce that $P_2(t_*)$ cannot be arbitrarily small. As before, using the fact that $P_2' \ge -CP_2$ it is seen that 
$P_2(t)\ge e^{-CT}P_2(t_*)$ and the conclusion follows. This, in turn, yields a lower bound $\mathfrak p_N>0$ for $P_N$. 

We are already in conditions of defining the open set 
$\Omega\subset \mathcal K^\circ$ as
$$\Omega:=\{ (M,P_0,P_1,P_2,P_N)\in C_T: \mathfrak m < M(t) < \mathcal M, 
\mathfrak{p}_0 < P_0(t) < \mathcal{P}_0,$$
$$\mathfrak{p}_1< P_1(t) < \mathcal{P}_1 ,\mathfrak{p}_2< P_2(t) < \mathcal{P}_2, \mathfrak{p}_N <P_N(t) < \mathcal{P}_N
\}$$
and 
\begin{theorem}
 \label{main}
Assume that the previous conditions (\ref{hip1}), 
(\ref{hip2}), (\ref{hip3}) and (\ref{hip4}) hold. 
Then problem (\ref{t1}) has at least 
one positive $T-$periodic solution.
 
\end{theorem}

\subsection{Degree computation}

In the previous section, the first condition of the continuation theorem was verified. It remains to prove that 
$b)$ and $c)$ are fulfilled as well. 
With this aim, set $\mathcal{Q}:=\Omega \cap \mathbb{R}^{5}$ and 
recall that the function 
$\phi:\overline {\mathcal{Q}}\to \mathbb{R}^5$ 
is defined by $\phi(x)=\overline {Nx}$. We claim that each coordinate 
$\phi_j$ has different signs at the corresponding opposite faces of $\mathcal Q$. 

Indeed, 
compute for example
$\phi_1(\mathcal{M},P_0,P_1,P_2,P_N)$ and $\phi_1(\mathfrak{m},P_0,P_1,P_2,P_N)$ for 
$\mathfrak{p}_j\le P_j \le \mathcal{P}_j$:
$$\phi_1(\mathcal{M},P_0,P_1,P_2,P_N)
= \frac 1T \int_0^T 
\left(\dfrac{V_S({t})K_1({t})^{n}}{K_1^{n}({t})+P_N}-
\dfrac{V_{m}({t})\mathcal{M}}{K_{m_1}({t})+\mathcal{M}}\right)\,dt$$
$$
<\mathcal V_S-\dfrac{v_{m}\mathcal{M}}{\mathcal K_{m_1}+\mathcal{M}}=0,$$
$$\phi_1(\mathfrak{m},P_0,P_1,P_2,P_N)
= \frac 1T \int_0^T 
\left(\dfrac{V_S({t})K_1({t})^{n}}{K_1^{n}({t})+P_N}-
\dfrac{V_{m}({t})\mathfrak{m}}{K_{m_1}({t})+\mathfrak{m}}\right)\,dt
$$
$$>\dfrac{v_S \kappa_1^{n}}{\mathcal K_1^{n}+\mathcal P_N^n}-
\dfrac{\mathcal V_{m}\mathfrak{m}}{\kappa_{m_1}+\mathfrak{m}}\ge 0$$
provided that $\mathfrak m$ is small enough.
In the same way, making the lower bounds smaller if necessary, 
we deduce that 
$$\phi_2(M,\mathcal P_0,P_1,P_2,P_N) 
< 0 < \phi_2(M,\mathfrak{p}_0,P_1,P_2,P_N)$$
$$\phi_3(M,P_0,\mathcal P_1,P_2,P_N) 
< 0 < \phi_2(M,{p}_0,\mathfrak{p}_1,P_2,P_N)$$
$$\phi_4(M,P_0,P_1,\mathcal P_2,P_N) 
< 0 < \phi_2(M,{p}_0,p_1,\mathfrak{p}_2,P_N)$$
$$\phi_5(M,P_0,P_1,P_2,\mathcal P_N) 
< 0 < \phi_2(M,{p}_0,p_1,p_2,\mathfrak{p}_N).$$
Thus, condition $b)$ of Continuation Theorem is verified. Moreover, we may define a homotopy as follows. 
Consider the center of  
$\mathcal{Q}$ given by
$$\wp:=\left(\dfrac{\mathcal{M}+\mathfrak{m}}{2},\dfrac{\mathcal{P}_0+\mathfrak{p}_0}{2},\dfrac{\mathcal{P}_1+\mathfrak{p}_1}{2},\dfrac{\mathcal{P}_2+\mathfrak{p}_2}{2},\dfrac{\mathcal{P}_N+\mathfrak{p}_N}{2} \right) $$
and the function 
$\mathcal{H}:\overline{\mathcal{Q}}\times [0;1]\to\R^5$ given by
$$\mathcal{H}(x,\lambda)=(1-\lambda)(\wp-x)+\lambda \phi.$$
We need to verify that $\mathcal H$ does not vanish at 
$\partial\mathcal Q$. To this end, suppose for example that 
$\mathcal{H}(\mathcal{M},P_0,P_1,P_2,P_N)=0$ for some $\hat{\lambda}\in[0;1]$, then
$$0=\mathcal{H}_{1}(\mathcal{M},\hat{\lambda})=(1-\hat{\lambda})
\underbrace{\left(\dfrac{\mathcal{M}+\mathfrak{m}}{2}-\mathcal{M}\right)}_{<0}+\hat{\lambda}\underbrace{\phi_1(\mathcal{M},P_0,P_1,P_2,P_N)}_{<0}<0,$$
a contradiction. All the remaining cases follow in an analogous way. 
By the homotopy invariance of the Brouwer degree, it follows that
$$deg_B(\phi,\mathcal{Q},0)= 
deg_B(\wp - I,\mathcal{Q},0)
=(-1)^{5}\neq 0.$$ 
This proves the third condition of the continuation theorem and, therefore, 
the existence of a $T$-periodic solution is deduced.\hfill$\square$\\

\end{document}